\documentclass[12pt]{article}
\newcommand{\al}{\alpha}
\newcommand{\be}{\beta}
\newcommand{\ga}{\gamma}
\newcommand{\g}[2]{g_{#1 \overline{#2}}}

\newcommand{\ddt}[1]{\frac{\partial #1}{\partial t}}
\newcommand{\dds}[1]{\frac{\partial #1}{\partial s}}

\newcommand{\F}{\mathcal{F}}
\newcommand{\imag}{\sqrt{-1}}
\newcommand{\subtitle}[1]{\bigskip \noindent
{\bf #1} \bigskip}
\begin{document}
\title{Higher K-energy functionals and higher Futaki invariants}
\newcounter{theor}
\setcounter{theor}{1}
\newtheorem{theorem}{Theorem}[theor]
\newtheorem{lemma}{Lemma}[section]
\newtheorem{remark}{Remark}[section]
\newtheorem{corollary}{Corollary}[section]
\newenvironment{proof}[1][Proof]{\begin{trivlist}
\item[\hskip \labelsep {\bfseries #1}]}{\end{trivlist}}

\centerline{\bf  Higher K-energy functionals and higher Futaki invariants}
\bigskip
\centerline{Ben Weinkove}
\bigskip
\centerline{Department of Mathematics}
\centerline{Columbia University, New York, NY 10027}
\bigskip
\bigskip

\subtitle{1. Introduction}

In 1986, Mabuchi \cite{ma} introduced the K-energy functional, which 
integrates the Futaki invariant \cite{fu} and whose critical points are metrics
of constant scalar curvature.  It was used by Bando and Mabuchi
\cite{bm} to prove a uniqueness theorem for K\"ahler-Einstein metrics.
The K-energy is strongly related to notions of stability in
Geometric Invariant Theory and has been used by Tian \cite{ti94,
ti97,ti} and Phong and Sturm \cite{ps}, to give some results on
the conjecture of Yau \cite{ya93,ya96} on the relationship between
stability and the existence of K\"ahler-Einstein metrics. Tian
\cite{ti00,ti} has cast both the K-energy and the Futaki invariant
in a more general setting using Bott-Chern forms and Donaldson
functionals.

Higher K-energy functionals were defined by Bando and Mabuchi 
\cite{bm86} and generalize 
the K-energy map to higher Chern classes of the manifold.   They 
integrate higher Futaki invariants (see \cite{bm86}, \cite{ba}, 
\cite{fu88}, \cite{fms}).  This note presents these functionals
with an emphasis on Bott-Chern forms.   Two new formulas for the higher 
K-energy 
 are given and the second K-energy is shown to be 
related to Donaldson's Lagrangian applied to metrics on the tangent 
bundle.

Let $M$ be a compact complex manifold of complex dimension $n$
with K\"ahler metric $\g{i}{j}$.  Let $\omega = \imag \g{i}{j} \,
dz^i \wedge dz^{\overline{j}}$ be the corresponding K\"ahler form.
Define
$$P(M, \omega) = \{\phi \in C^{\infty}(M) \, | \, \g{i}{j} +
 \partial_i \partial_{\overline{j}} \phi >0 \}.$$  For $\phi$
in $P(M, \omega)$, let $\{ \phi_t \}_{0 \le t \le 1}$ be a smooth
path in $P(M, \omega)$ with $\phi_0 = 0$, $\phi_1 = \phi$.  Write
$\omega_t$ for $\omega + \imag \partial \overline{\partial} \phi_t$.
Now define the \emph{$k$th $K$-energy functional} by
$$ M_{k, \omega} (\phi) = - (n-k+1) \int_0^1 \int_M \dot{\phi_t} (c_k
(\phi_t) - \mu_k \omega_t^k) \wedge \omega_t^{n-k} dt,$$ where
$c_k (\phi_t)$ (also written $c_k(\omega_t)$)  is the representative of the 
$k$th Chern class of
$M$ given by the metric $\omega_t$, and $\mu_k$ is the number given by
$$\mu_k = \frac{1}{V} \int_M c_k (\omega) \wedge \omega^{n-k}, \quad
\textrm{where} \quad V= \int_M \omega^n.$$  The case $k=1$
corresponds to the original K-energy as defined in \cite{ma}.

Theorems 1 and 2 below were proved by Bando and Mabuchi \cite{bm86} 
(also see \cite{fu88} for a proof of Theorem 1).  We present an 
alternative 
derivation using Bott-Chern 
forms, leading to new formulas for the higher K-energy (Lemma 5.1 and 
Corollary 5.2) as well as Theorem 3 below, which gives the second K-energy 
in 
terms of Donaldson's Lagrangian.

\bigskip \noindent
{\bf Theorem 1.}  \ \ 
$M_{k, \omega} (\phi)$ is independent of the
choice of path $\{\phi_t\}$. \bigskip

To define the higher Futaki invariants, let $c_k (\omega)$ and
$Hc_k (\omega)$ be respectively the $k$th Chern form and its
harmonic part.  There is a real smooth $(k-1, k-1)$ form $f_k$
such that
$$c_k (\omega) - H c_k (\omega) = \imag \partial
\overline{\partial} f_k.$$ Define the \emph{kth Futaki invariant},
which acts on a holomorphic vector field $X$, by
$$\F_{k, \omega} (X) = \int_M \mathcal{L}_X f_k \wedge
\omega^{n-k+1},$$ where $\mathcal{L}_X$ is the Lie derivative. The
case $k=1$ corresponds to the original Futaki invariant if $c_1(M)
= \mu_1 \omega$, so that $Hc_1(\omega)=\mu_1 \omega$. It is
shown in \cite{ba} (see also \cite{fu88}) that
$\F_{k,\omega}$ is independent of choice of metric in the K\"ahler
class.

Now given a holomorphic vector field $X$, let $\Phi_t$ be the
integral curve of $X_R = X + \overline{X}$.  Then $\Phi_t : 
M
\rightarrow M$ is holomorphic, and there exists a smooth path $\{
\phi_t \}$ in $P(M, \omega)$ satisfying
$$\Phi_t^* \omega = \omega + \imag \partial \overline{\partial}
\phi_t, \quad \textrm{and} \quad \int_M \phi_t \omega^n =0.$$
Then the higher K-energy integrates the Futaki invariant in the following 
sense. 

\bigskip \noindent {\bf Theorem 2.} \ \
Suppose that $c_k (M) = \mu_k [ \omega^k ]$. Then with $\{\phi_t \}$
as above,
$$ \frac{d}{dt} M_{k, \omega} (\phi_t) = 2 \textrm{Re} ( 
\F_{k, \omega} (X)).$$

Donaldson \cite{do,do2} defined a Lagrangian which was used in his proof that 
Mumford-Takemoto stable vector bundles over projective algebraic varieties 
admit Hermitian-Einstein metrics. 
It is defined in terms of the Bott-Chern forms corresponding to $\textrm{Tr}(A)
$ and $\textrm{Tr}(AB)$.  It is unsurprising then that this Lagrangian is 
related to the first and second K-energy functionals, which correspond to the 
first and second Chern classes.  

Let $E$ be a holomorphic vector bundle of rank $r$ over $M$, and let $H$ and 
$H_0$ be Hermitian metrics on $E$.  Let $\{H_t\}_{0 \le t \le 1 }$ be a smooth 
path of metrics between $H_0$ and $H_1 = H$.  Let $F_t$ be the curvature of 
$H_t$.  Then Donaldson's Lagrangian is given by
\begin{eqnarray*}
\lefteqn{L(H, H_0) } \\
&& \mbox{} = \int_M \int_0^1 n \imag \textrm{Tr} (\dot{H}_t H_t^{-1} F_t) 
\wedge 
\omega^{n-1}dt - \lambda \int_M \int_0^1 \textrm{Tr} (\dot{H}_t H_t^{-1}) 
\omega^n dt,
\end{eqnarray*}
where
$$ \lambda = \frac{2\pi n}{r} \frac{1}{V}\int_M c_1(E) \wedge \omega^{n-1}.$$
Now consider the case when $E = T'M$, the holomorphic tangent bundle.  Let $H_0
$ be the K\"ahler metric $g$, and, for $\phi$ in $P(M,\omega)$, let $H$ be $g_
{\phi}$.  Write $\omega_{\phi}$ for $\omega + \imag \partial \overline
{\partial} \phi$ and $L(\phi)$ for $L(g_\phi, g)$.  We have the 
following formula for the second K-energy functional.

\bigskip \noindent {\bf Theorem 3.} \ \ For $\phi$ in $P(M,\omega)$,
\begin{eqnarray*}
\lefteqn{M_{2, \omega} (\phi) = - \frac{1}{(2 \pi)^2 n} L(\phi) }\\
&& \mbox{} - \frac{1}{4 \pi n} \int_M \log \left( \frac{\omega_{\phi}^n}
{\omega^n} \right) (2 \mu_1 \omega - n c_1(\phi) - n c_1 (0))\wedge \omega^{n-
1} \\
&& \mbox{} - \int_M \phi \left( c_2 (\phi) \wedge \sum_{i=0}^{n-2} \omega^i 
\wedge \omega_{\phi}^{n-2-i} - \left( \frac{n-1}{n+2} \right) \mu_2 \sum_{i=0}
^n \omega^i \wedge \omega_{\phi}^{n-i} \right).
\end{eqnarray*}

I am very grateful to my advisor, D.H. Phong, for suggesting to look at 
these 
functionals, and for his subsequent guidance.  I would also like to thank all 
the members of the Geometry and Analysis community at Columbia University for 
many enlightening discussions and talks.  This paper will form a part of my 
forthcoming PhD thesis at Columbia University.

\bigskip
\bigskip
\pagebreak
\subtitle{2. Chern classes and Bott-Chern forms}

Before proving the above theorems, we will briefly discuss Chern classes 
and 
their 
Bott-Chern forms.  For more details on Bott-Chern forms, see \cite{bc}, \cite
{do} and
\cite{ti}. Let $E$ be a holomorphic vector bundle of rank $r$ over
the compact K\"ahler manifold $M$. If $H=H_{\al \overline{\be}}$
is a Hermitian metric on $E$, then its curvature is an
endomorphism-valued (1,1) form given by
$$(F_H)_{\al}^{\ \be} = {(F_H)_{\al}^{\ \be}}_{i \overline{j}}\
dz^i \wedge dz^{\overline{j}} = \overline{\partial} ((\partial H_{\al \overline
{\ga}}) H^{\overline{\ga} \be}).$$ Let $\Phi$ be an invariant
symmetric k-linear function on $\mathbf{gl}(r, \mathbf{C})$.  Then
we have a Chern-Weil form
$$\Phi(H) = \Phi(F_H)= \Phi(F_H , \ldots ,
F_H) \in \Omega_M^{k,k},$$ which represents a characteristic class.

Let $H$ and $H_0$ be two Hermitian metrics on $E$, and let
$\{H_t\}$ be a smooth path of metrics between them. Then there exists a Bott-
Chern form $\textrm{BC}_{\Phi}(H, H_0)$ in
$\Omega_M^{k-1,k-1} / (\textrm{Im}\,
\partial + \textrm{Im}\, \overline{\partial})$ given by
$$\textrm{BC}_{\Phi}(H,H_0) = - k\imag \int_0^1 \Phi (\dot{H}_t H_t^{-1},
F_{H_t}, \ldots, F_{H_t}) dt.$$ It is shown in \cite{bc} that this
definition is independent of the choice of path.  Moreover we have
$$ - \imag \partial \overline{\partial} \textrm{BC}_{\Phi}(H, H_0) =
\Phi(H) - \Phi(H_0),$$ which can be checked by differentiating
with respect to some parameter.

We will restrict to the case of Chern classes of holomorphic
vector bundles.  The Chern-Weil form for the metric $H$
representing the $k$th Chern class of $E$ is given by
$$c_k (H) = \frac{1}{k!}\left( \frac{\imag}{2\pi}\right)^k
\sum_{\pi \in \Sigma_k} \textrm{sgn}(\pi)\, F_{\al_{\pi(1)}}^{\
\al_1} \wedge F_{\al_{\pi(2)}}^{\ \al_2} \wedge \ldots \wedge
F_{\al_{\pi(k)}}^{\ \al_k},$$ where $F=F_H$ and the usual
summation convention is being used.

The Bott-Chern form for the $k$th Chern class is then given by
\begin{eqnarray*}
\lefteqn{\textrm{BC}_k(H, H_0)} \\
&& = -k\imag \int_0^1 \frac{1}{k!}\left(
\frac{\imag}{2\pi}\right)^k \sum_{\pi \in \Sigma_k}
\textrm{sgn}(\pi)\, (\dot{H_t} H_t^{-1})_{\al_{\pi(1)}}^{\ \al_1}
F_{\al_{\pi(2)}}^{\ \al_2} \wedge \ldots \wedge
F_{\al_{\pi(k)}}^{\ \al_k}dt,
\end{eqnarray*}
where $F=F_{H_t}.$  It satisfies
$$ c_k(H) - c_k(H_0) = - \imag \partial \overline{\partial} \textrm{BC}_k
(H,H_0).$$

Now consider the case where the vector bundle is $T'M$, the
holomorphic tangent bundle of $M$.  Then $\g{i}{j}$, the K\"ahler
metric on $M$, is a Hermitian metric on $T'M$.  Let $\phi$ be in
$P(M, \omega)$.  Write $g_{\phi}$ for $\g{i}{j} + \partial_i
\partial_{\overline{j}}
\phi$ and $c_k (\phi)$ for $c_k (g_\phi)$.  Then, in normal
coordinates for $g_\phi$, we have
$$c_k (\phi) = \frac{1}{k!}\left( \frac{\imag}{2\pi}\right)^k \sum_{\pi \in 
\Sigma_k}
\textrm{sgn}(\pi)\, F_{\al_{\pi(1)} \overline{\al}_1} \wedge
 \ldots \wedge F_{\al_{\pi(k)} \overline{\al}_k},$$
writing $F_{\al \overline{\be}}$ for $R_{\al \overline{\be} i
\overline{j}} dz^i \wedge dz^{\overline{j}}$, where $R_{\al
\overline{\be} i \overline{j}}$ is the curvature tensor for
$g_{\phi}$.

\setcounter{section}{3}
\setcounter{lemma}{0}
\bigskip
\bigskip \noindent
{\bf 3. Higher K-energy functionals}
\bigskip

In this section we give an alternative proof of Theorem 1.
Let $\{ \phi_t \}_{0\le t \le 1}$ be a smooth path in
$P(M,\omega)$ with $\phi_0 = \phi_1$, and let $\lambda$ be any
constant.  Calculate
\begin{eqnarray*}
\lefteqn{\int_0^1 \int_M \dot{\phi_t} (c_k (\phi_t) - \lambda
\omega_t^k) \wedge \omega_t^{n-k} dt} \\
& = & \int_0^1 \int_M \dot{\phi_t} (c_k(\phi_t) - c_k(\phi_0) +
c_k(\phi_0) - \lambda \omega_t^k)\wedge \omega_t^{n-k} dt \\
& = & \int_0^1 \int_M \dot{\phi_t} (- \imag \partial
\overline{\partial} \textrm{BC}_k (\phi_t,\phi_0) + c_k(\phi_0) -
\lambda \omega_t^k)
\wedge \omega_t^{n-k} dt \\
& = & \int_0^1 \int_M (- \imag \partial \overline{\partial}
\dot{\phi_t}) \wedge \textrm{BC}_k (\phi_t, \phi_0) \wedge
\omega_t^{n-k}
dt \\
&& \mbox{} + \int_0^1 \int_M \dot{\phi_t} (c_k(\phi_0) - \lambda
\omega_t^k) \wedge \omega_t^{n-k} dt\\
& = & - \frac{1}{n-k+1}\int_0^1 \int_M  \textrm{BC}_k(\phi_t,
\phi_0) \wedge \ddt{} (\omega_t^{n-k+1}) dt \\
&& \mbox{} + \int_0^1 \int_M \dot{\phi_t} (c_k(\phi_0) - \lambda
\omega_t^k)\wedge \omega_t^{n-k} dt\\
& = & \frac{1}{n-k+1} \int_0^1 \int_M \ddt{} \textrm{BC}_k
(\phi_t, \phi_0) \wedge \omega_t^{n-k+1} dt \\
&& \mbox{} + \int_0^1 \int_M \dot{\phi_t} (c_k (\phi_0) - \lambda
\omega_t^k) \wedge \omega_t^{n-k} dt.
\end{eqnarray*}
The next two lemmas will show that this expression is zero.  The
proof of the theorem will then follow immediately.

\begin{lemma} \label{lemmawelldefined1} With $\{ \phi_t \}$ as above,
$$\int_0^1 \int_M \dot{\phi_t} (c_k(\phi_0) - \lambda \omega_t^k)
\wedge \omega_t^{n-k} dt = 0.$$
\end{lemma}
\begin{proof}
The argument is almost identical to one given in \cite{ti94}, but we
will include it for the reader's convenience.  Define
$$\phi_{t,s} = s \phi_t, \quad \textrm{and} \quad \omega_{t,s} = \omega +
\imag\partial \overline{\partial} \phi_{t,s},$$
$$ D(s) = \int_0^1 \int_M \ddt{\phi_{t,s}} (c_k(\phi_0) - \lambda
\omega_{t,s}^k) \wedge \omega_{t,s}^{n-k} dt.$$ Notice that
$D(0)=0$.  We will show that $D'(s)=0$.
\setlength\arraycolsep{2pt}
\begin{eqnarray*}
D'(s) & = & \int_0 \int_M \frac{\partial^2 \phi_{t,s}}{\partial t
\partial s} (c_k (\phi_0) - \lambda \omega_{t,s}^k)
\wedge \omega_{t,s}^{n-k}dt \\
&& \mbox{} + \int_0^1 \int_M \ddt{\phi_{t,s}} ((n-k) c_k(\phi_0)
\omega_{t,s}^{n-k-1} - \lambda n \omega_{t,s}^{n-1}) \wedge \imag
\partial \overline{\partial} (\ddt{\phi_{t,s}}) dt\\
& = & \int_0^1 \int_M \ddt{} (\dds{\phi_{t,s}}) (c_k (\phi_0) -
\lambda \omega_{t,s}^k) \wedge \omega_{t,s}^{n-k} dt\\
&& \mbox{} + \int_0^1 \int_M \dds{\phi_{t,s}} ((n-k) c_k(\phi_0)
\omega_{t,s}^{n-k-1} - \lambda n \omega_{t,s}^{n-1})\wedge \imag
\partial \overline{\partial} (\ddt{\phi_{t,s}}) dt\\
& = & \int_0^1 \int_M \ddt{} \left( \dds{\phi_{t,s}} (c_k(\phi_0)
- \lambda \omega_{t,s}^k) \wedge \omega_{t,s}^{n-k} \right) dt \\
& = & 0.
\end{eqnarray*}
\end{proof}
\begin{lemma} \label{lemmaBC} With $\{\phi_t \}$ as above,
$$\int_M \ddt{} ( \emph{BC}_k(\phi_t , \phi_0) ) \wedge \omega_t^{n-k+1}
=0.$$
\end{lemma}
\begin{proof}
Writing $g$ for $g_{\phi_t}$ we have
$$\dot{g}_{i \overline{j}} g^{\overline{j} k} = (\partial_i
\partial_{\overline{j}} \dot{\phi_t}) g^{\overline{j} k}.$$
Hence, using the formula for the Bott-Chern form, and working in
normal coordinates, \setlength\arraycolsep{2pt}
\begin{eqnarray}
\nonumber \lefteqn{\int_M \ddt{} \textrm{BC}_k(\phi_t, \phi_0)
\wedge\omega_t^{n-k+1}}\\
\nonumber & = & \int_M \frac{- \imag}{(k-1)!} \left(
\frac{\imag}{2\pi} \right)^k \sum_{\pi \in \Sigma_k}
\textrm{sgn}(\pi)
\partial_{\al_{\pi(1)}} \partial_{\overline{\al}_1} \dot{\phi_t}
F_{\al_{\pi(2)} \overline{\al}_2} \wedge \ldots \\
\nonumber && \qquad \qquad \qquad \mbox{} \ldots \wedge
F_{\al_{\pi(k)} \overline{\al}_k}\wedge \omega_t^{n-k+1}\\
\nonumber & = & \frac{-\imag}{(k-1)!}\left( \frac{\imag}{2\pi}
\right)^k \int_M \sum_{\pi \in \Sigma_k} \textrm{sgn}(\pi)
\nabla_{\al_{\pi(1)}} \nabla_{\overline{\al}_1} \dot{\phi_t}
F_{\al_{\pi(2)}
\overline{\al}_1} \wedge \ldots \\
&& \qquad \qquad \qquad \mbox{} \ldots
 \wedge F_{\al_{\pi(k)}
\overline{\al}_k}\wedge \omega_t^{n-k+1}. \label{integral}
\end{eqnarray}
Integrating by parts, using the Bianchi identity, and defining
$\tau$ to be the transposition in $\Sigma_k$ which interchanges 1
and 2, we see that one term in the integrand of (\ref{integral}) will be
\begin{eqnarray*}
\lefteqn{- \sum_{\pi \in \Sigma_k} \textrm{sgn} (\pi)
\nabla_{\overline{\al}_1} \dot{\phi_t} \nabla_{\al_{\pi (1)}}
F_{\al_{\pi(2)} \overline{\al}_2} \wedge F_{\al_{\pi(3)}
\overline{\al}_3} \wedge \ldots \wedge F_{\al_{\pi(k)}
\overline{\al}_k} \wedge
\omega_t^{n-k+1}} \\
&= & - \sum_{\pi \in \Sigma_k} \textrm{sgn} (\pi)
\nabla_{\overline{\al}_1} \dot{\phi_t}
 \nabla_{\al_{\pi (2)}} F_{\al_{\pi(1)}
\overline{\al}_2} \wedge F_{\al_{\pi(3)} \overline{\al}_3} \wedge
\ldots \wedge F_{\al_{\pi(k)} \overline{\al}_k} \wedge
\omega_t^{n-k+1}\\
& = & - \sum_{\pi \in \Sigma_k} \textrm{sgn} (\pi \tau)
\nabla_{\overline{\al}_1} \dot{\phi_t}  \nabla_{\al_{\pi (1)}}
F_{\al_{\pi(2)} \overline{\al}_2} \wedge F_{\al_{\pi(3)}
\overline{\al}_3} \wedge \ldots \wedge F_{\al_{\pi(k)}
\overline{\al}_k} \wedge
\omega_t^{n-k+1}\\
& = & \sum_{\pi \in \Sigma_k} \textrm{sgn} (\pi)
\nabla_{\overline{\al}_1} \dot{\phi_t}
 \nabla_{\al_{\pi (1)}} F_{\al_{\pi(2)}
\overline{\al}_2} \wedge F_{\al_{\pi(3)} \overline{\al}_3} \wedge
\ldots \wedge F_{\al_{\pi(k)} \overline{\al}_k} \wedge
\omega_t^{n-k+1}\\
& = & 0.
\end{eqnarray*}
Every other term is zero by the same argument, and the lemma is
proved.
\end{proof}

The proof of Theorem 1 is now complete.

\bigskip
\noindent
{\bf Remark 3.3} \  $M_{k,\omega}$ satisfies a cocycle condition.  Namely, 
if $\omega'$ is another K\"ahler metric with $\omega' = \omega + \imag 
\partial \overline{\partial} \psi$, then
$$M_{k, \omega} (\phi) - M_{k,\omega'} (\phi - \psi) = M_{k, \omega} 
(\psi).$$

\bigskip
\noindent
{\bf Remark 3.4} \
The critical points of $M_{k, \omega}$ are metrics $\omega$ with 
$\Lambda^k c_k (\omega)$ constant.  An example of such a metric is the 
Fubini-Study metric on $\mathbf{CP}^n$.

\setcounter{section}{4}
\setcounter{lemma}{0}
\bigskip
\bigskip \noindent
{\bf 4. Higher Futaki invariants} \bigskip

In this section we give a proof of Theorem 2.
Let $X$ be a holomorphic vector field.  Recall that the Lie
derivative of $X$ on forms is given by
$$\mathcal{L}_X = i_X \circ d + d \circ i_X.$$
The interior product $i_X \omega$ is a (0,1) form which is
$\overline{\partial}$-closed since $\omega$ is K\"ahler. Hence
there exists a smooth function $\theta_X$ and a
harmonic $(0,1)$ form $\al$ such that
$$i_X \omega = \al - \overline{\partial} \theta_X.$$
We have the following lemma.

\begin{lemma}
$$\F_{k, \omega } (X) = -(n-k+1) \imag \int_M \theta_X (c_k (\omega) -
Hc_k(\omega) ) \wedge \omega^{n-k},$$ where $\theta_X$ is as given
above.
\end{lemma}
\begin{proof}  Integrating by parts,
\setlength\arraycolsep{2pt}
\begin{eqnarray*}
\F_{k,\omega} (X) & = & \int_M \mathcal{L}_X f_k \wedge \omega^{n-k+1} \\
& = & - \int_M f_k \wedge \mathcal{L}_X \omega^{n-k+1} \\
& = & - \int_M f_k \wedge \partial i_X \omega^{n-k+1} \\
& = & - (n-k+1) \int_M f_k \wedge \partial (\al - 
\overline{\partial}\theta_X 
) \wedge \omega^{n-k} \\ 
& = & (n-k+1) \int_M f_k \wedge \partial \overline{\partial} \theta_X 
\wedge \omega^{n-k} \\
& = & -(n-k+1) \imag \int_M \theta_X (c_k(\omega) - H c_k(\omega)) \wedge 
\omega^{n-k}.
\end{eqnarray*}
\end{proof}

To finish the proof of the theorem, let $\Phi_t$ be the integral
curve of $X_R$ as before, and let $\phi_t$ be given
by
$$\Phi_t^* \omega = \omega_t = \omega + \imag \partial \overline{\partial}
\phi_t, \qquad \int_M \phi_t \omega^n = 0. $$ Then
$$\mathcal{L}_{X_R} \omega_t = \imag \partial 
\overline{\partial} \dot{\phi_t}
.$$ Hence
$$ \partial i_X \omega_t + \overline{\partial} i_{\overline{X}} \omega_t = 
\frac{1}{2} \imag (\partial \overline{\partial} \dot{\phi}_t - 
\overline{\partial} \partial \dot{\phi}_t)$$
and so
$$i_X \omega_t - \frac{1}{2} \imag \overline{\partial} \dot{\phi}_t =
\al_t,$$ where $\al_t$ is a (0,1) form satisfying $\textrm{Re} (\partial 
\al_t)=0$.  It follows from above that
\setlength\arraycolsep{2pt}
\begin{eqnarray*}
2 \textrm{Re} (\F_{k, \omega_t} (X)) & = & - (n-k+1) \int_M
\dot{\phi_t}
(c_k (\omega_t) - Hc_k (\omega_t)) \wedge \omega_t^{n-k} \\
& = & - (n-k+1) \int_M \dot{\phi_t}
 (c_k
(\omega_t)
- \mu_k \omega_t^k) \wedge \omega_t^{n-k} \\
& = &  \frac{d}{dt} M_{k, \omega} (\phi_t).
\end{eqnarray*}
Finally we must show that $\F_{k, \omega_t}(X) = \F_{k,\omega}
(X)$. This is immediate from \cite{ba}, or alternatively can be seen as 
follows.  Since 
$$ \mathcal{L}_{X_R} \omega_t = \Phi_t^* 
\mathcal{L}_{X_R} \omega,$$
we have
$$\partial \overline{\partial} \dot{\phi}_t = \partial \overline{\partial} 
\Phi_t^* \dot{\phi}_0,$$
and so $\dot{\phi}_t = \Phi_t^* \dot{\phi}_0$.  Hence
\begin{eqnarray*}
\int_M \dot{\phi}_t (c_k(\omega_t) - \mu_k \omega_t^k) \wedge
\omega_t^{n-k} & = & \int_M \Phi_t^* \dot{\phi}_0 (\Phi_t^* c_k
(\omega) - \mu_k \Phi_t^* \omega^k)\wedge \Phi_t^* \omega^{n-k} \\
& = & \int_M \dot{\phi}_0 (c_k (\omega) - \mu_k \omega^k) \wedge
\omega^{n-k}.
\end{eqnarray*}

\setcounter{section}{5}
\setcounter{lemma}{0}
\subtitle{5. Higher K-energy functionals and Donaldson's Lagrangian}

In this section we give a proof of Theorem 3 and give two 
alternative formulas for the higher K-energy functionals.
It is pointed out in \cite{ti00} (see also \cite{ch} and \cite{ti94}) that the 
K-energy can be written without a path integral as
\begin{eqnarray*}
\lefteqn{M_{1,\omega} (\phi) = \int_M \frac{1}{2\pi} \log \left( \frac{\omega_
{\phi}^n}{\omega^n} \right) \omega_{\phi}^n}\\ 
&& \mbox{} - \int_M \phi \left( \frac{\imag}{2\pi} \textrm{Ric}(\omega) \wedge 
\sum_{i=0}^{n-1} \omega^i \wedge \omega_{\phi}^{n-1-i}  - \frac{\mu_1}{n+1} 
\sum_{i=0}^n \int_M \phi \omega^i \wedge \omega_{\phi}^{n-i}\right).
\end{eqnarray*}
This appears to be both a natural and useful expression for the K-energy (see 
\cite{ps}).
Notice that
$$ c_1(\omega) = \frac{\imag}{2\pi} \textrm{Ric}(\omega) \quad \textrm{and} 
\quad
\frac{1}{2\pi} \log  \left( \frac{\omega_{\phi}^n}{\omega^n} \right) = \textrm
{BC}_1 (\phi,0),$$
where to be precise, the above is only a \emph{representative} of $\textrm{BC}
_1 (\phi, 0)$. 
The following lemma generalizes the above formula to the higher K-energy 
functionals.  

\begin{lemma} \label{lemmaformula}
\begin{eqnarray*}
\lefteqn{M_{k, \omega} (\phi) = \int_M \textrm{BC}_k (\phi, 0) \wedge 
\omega_\phi^{n-k+1} } \\
&& \mbox{} - \int_M \phi \left(c_k(0) \wedge \sum_{i=0}^{n-k} \omega^i \wedge 
\omega_{\phi}^{n-k-i} - \frac{(n-k+1)}{n+1} \mu_k \sum_{i=0}^n \omega^i \wedge 
\omega_{\phi}^{n-i} \right).
\end{eqnarray*}
\end{lemma}

\begin{proof} For a path $\{ \phi_t \}$ in $P(M, \omega)$, calculate
\setlength\arraycolsep{2pt}
\begin{eqnarray}
\nonumber  \lefteqn{\int_M \ddt{} ( \textrm{BC}_k (\phi_t, 0) \wedge \omega_t^
{n-k+1})}\\
\nonumber
&& \qquad \qquad =  \int_M \ddt{} (\textrm{BC}_k (\phi_t, 0)) \wedge \omega_t^
{n-k+1} \\ 
\nonumber
&& \qquad \qquad \mbox{} + (n-k+1) \int_M \textrm{BC}_k (\phi_t,0) \wedge 
\omega_t^{n-k} \wedge \imag 
\partial \overline{\partial} \dot{\phi}_t \\
&& \qquad \qquad =  - (n-k+1) \int_M \dot{\phi}_t (c_k (\phi_t) - c_k(0)) 
\wedge \omega_t^{n-k},
\label{part1}
\end{eqnarray}
using Lemma \ref{lemmaBC} for the first term and integrating by parts for the 
second.  Also
\begin{eqnarray}
\nonumber
\lefteqn{\int_M \ddt{} ( \phi_t c_k(0) \wedge \sum_{i=0}^{n-k} \omega^i \wedge 
\omega_t^{n-k-i} )} \\
\nonumber
&& \qquad = \int_M \dot{\phi}_t c_k(0) \wedge \sum_{i=0}^{n-k} \omega^i \wedge 
\omega_t^{n-k-i} \\
\nonumber
&&\qquad \mbox{} + \int_M \phi_t c_k(0) \wedge \sum_{i=0}^{n-k} (n-k-i) 
\omega^i \wedge  \omega_t^{n-k-i-1} \wedge \imag \partial \overline{\partial} 
\dot{\phi}_t \\
\nonumber
&& \qquad  = \int_M \dot{\phi}_t c_k(0) \wedge \sum_{i=0}^{n-k} \omega^i \wedge 
\omega_t^{n-k-i} \\
\nonumber
&& \qquad \mbox{} + \int_M \dot{\phi}_t c_k(0) \wedge \sum_{i=0}^{n-k} (n-k-i) 
\omega^i \wedge \omega_t^{n-k-i-1} \wedge (\omega_t - \omega) \\
\nonumber
&& \qquad = \int_M \dot{\phi}_t c_k(0) \wedge ( \sum_{i=0}^{n-k} (n-k-i+1) 
\omega^i \wedge \omega_t^{n-k-i} \\
\nonumber
&& \qquad \mbox{} - (n-k-i) \omega^{i+1} \wedge \omega_t^{n-k-i-1} )\\
&& \qquad = (n-k+1) \int_M \dot{\phi}_t c_k(0) \wedge \omega_t^{n-k}. \label
{part2}
\end{eqnarray}
Similarly,
\begin{eqnarray}
\int_M \ddt{} \left( \frac{(n-k+1)\mu_k}{n+1} \phi_t \sum_{i=0}^{n} \omega^i 
\wedge \omega_t^{n-i}\right) = (n-k+1) \mu_k \int_M \dot{\phi}_t \omega_t^n. 
\label{part3}
\end{eqnarray}
The lemma follows from (\ref{part1}), (\ref{part2}) and (\ref{part3}).
\end{proof}

We will need a slightly different formula for $M_{k,\omega}$.

\setcounter{corollary}{1}
\begin{corollary} \label{corollarykenergy}
\begin{eqnarray*}
\lefteqn{M_{k, \omega} (\phi) = \int_M \textrm{BC}_k (\phi, 0) \wedge \omega^{n-
k+1} } \\
&& \mbox{} - \int_M \phi \left(c_k(\phi) \wedge \sum_{i=0}^{n-k} \omega^i 
\wedge \omega_{\phi}^{n-k-i} - \frac{(n-k+1)}{n+1} \mu_k \sum_{i=0}^n \omega^i 
\wedge \omega_{\phi}^{n-i} \right).
\end{eqnarray*}
\end{corollary}
\addtocounter{lemma}{1}
\begin{proof}
\begin{eqnarray*}
\lefteqn{\int_M \textrm{BC}_k (\phi,0) \wedge \omega_{\phi}^{n-k+1}} \\
& = & \int_M \textrm{BC}_k (\phi,0) \wedge (\omega + \imag \partial \overline
{\partial}\phi)^{n-k+1} \\
& = & \int_M \textrm{BC}_k (\phi,0) \wedge \sum_{i=0}^{n-k+1} {n-k+1 \choose i} 
\omega^i \wedge  (\omega_{\phi} - \omega)^{n-k+1-i} \\
& = & \int_M \textrm{BC}_k (\phi ,0) \wedge \omega^{n-k+1} \\
&& \mbox{} + \int_M \textrm{BC}_k (\phi, 0) \wedge \imag \partial \overline
{\partial} \phi
\wedge \sum_{i=0}^{n-k} {n-k+1 \choose i} \omega^i \wedge (\omega_{\phi} - 
\omega)^{n-k-i} \\
& = & \int_M \textrm{BC}_k (\phi,0) \wedge \omega^{n-k+1} - \int_M \phi (c_k 
(\phi) - c_k(0)) \wedge \sum_{i=0}^{n-k} \omega^i \wedge \omega_{\phi}^{n-k-i},
\end{eqnarray*}
using, for the last line, the elementary identity
$$\sum_{i=0}^m {m+1 \choose i} x^i (y-x)^{m-i} = \sum_{i=0}^m x^i y^{m-i},$$
which can be seen by observing that each side is equal to 
$$(y^{m+1} - x^{m+1})/(y-x).$$
The corollary now follows immediately from Lemma \ref{lemmaformula}.
\end{proof}

We now prove Theorem 3.  First note that we have
\begin{eqnarray*}
\textrm{BC}_1 (\phi , 0) & = & \frac{1}{2\pi} \int_0^1 \textrm{Tr} (\dot{g}_t 
g_t^{-1})dt, \qquad \textrm{and} \\
\textrm{BC}_2 (\phi , 0) & = & \frac{\imag}{(2\pi)^2} \int_0^1 \left( \textrm
{Tr} (\dot{g}_t g_t^{-1}) \textrm{Tr}(F_t) - \textrm{Tr} (\dot{g}_t g_t^{-1} 
F_t) \right) dt,
\end{eqnarray*}
where $F_t$ is the curvature of $g_t$.  Then using the fact that $\lambda = 2
\pi \mu_1$ we have
\begin{eqnarray}
\nonumber
L(\phi) & = & n (2\pi)^2 \int_M \int_0^1 \frac{\imag}{(2\pi)^2} \textrm{Tr} 
(\dot{g}_t g_t^{-1} F_t) \wedge \omega^{n-1} dt\\
\nonumber
&& \mbox{} - (2\pi)^2 \mu_1 \int_M \int_0^1 \frac{1}{2\pi} \textrm{Tr} (\dot{g}
_t g_t^{-1}) \omega^n dt \\
\nonumber 
& = & -n (2\pi)^2 \int_M \textrm{BC}_2 (\phi, 0) \wedge \omega^{n-1} - (2\pi)^2 
\mu_1 \int_M \textrm{BC}_1 (\phi,0) \omega^n \\
&& \mbox{} + n\imag \int_M \int_0^1 \textrm{Tr} (\dot{g}_t g_t^{-1})\textrm{Tr} 
(F_t) \wedge \omega^{n-1} dt. \label{donaldson}
\end{eqnarray}
The last term is equal to
\begin{eqnarray}
\nonumber
\lefteqn{n\imag \int_M \int_0^1 (\triangle_t \dot{\phi}_t) \textrm{Ric} 
(\omega_t) \wedge \omega^{n-1} dt}\\
\nonumber && \mbox{} = \frac{n \imag}{2} \int_M \log \left( \frac{\omega_{\phi}
^n}{\omega^n} \right) (\textrm{Ric} (\omega_t) + \textrm{Ric} (\omega)) \wedge 
\omega^{n-1} \\
&& \mbox{} = n \pi \int_M \log \left( \frac{\omega_{\phi}^n}{\omega^n} \right) 
(c_1(\phi) + c_1(0) ) \wedge \omega^{n-1},  \label{lastterm}
\end{eqnarray}
where a straightforward integration by parts has been used for the second 
line.
The proof of the theorem now follows from Corollary \ref{corollarykenergy} and 
equations (\ref{donaldson}) and (\ref{lastterm}).

\end{document}